\documentclass[11pt,a4paper]{article}
\usepackage{amsfonts,amsmath,amssymb,accents,amsthm}
\usepackage{verbatim}

\usepackage{hyperref}
\usepackage[normalem]{ulem}

\usepackage{tikz} 

\newcommand{\dis}{\displaystyle}
\newcommand{\txt}{\textstyle}

\newcommand{\noi}{\noindent}
\newcommand{\halmos}{\rule{1ex}{1.4ex}}
\newcommand{\QED}{\nopagebreak{\hspace*{\fill}$\halmos$\medskip}}

\newcommand{\quand}{\quad\mbox{and}\quad}

\newtheoremstyle{mythm}
  {}
  {}
  {\itshape}
  {}
  {\bfseries}
  {}
  {.5em}
  {#1 #2 \thmnote{(#3)}}

\theoremstyle{mythm}
\newtheorem{theorem}{Theorem}
\newtheorem{proposition}[theorem]{Proposition}
\newtheorem{lemma}[theorem]{Lemma}
\newtheorem{exercise}[theorem]{Exercise}
\newtheorem{corollary}[theorem]{Corollary}
\newtheorem{conjecture}[theorem]{Conjecture}

\newtheorem{counterex}[theorem]{Counterexample}

\newcommand{\bt}{\begin{theorem}}
\newcommand{\et}{\end{theorem}}
\newcommand{\bl}{\begin{lemma}}
\newcommand{\el}{\end{lemma}}
\newcommand{\bp}{\begin{proposition}}
\newcommand{\ep}{\end{proposition}}
\newcommand{\bcor}{\begin{corollary}}
\newcommand{\ecor}{\end{corollary}}
\newcommand{\br}{\begin{remark}\rm}
\newcommand{\er}{\end{remark}}
\newcommand{\bcon}{\begin{conjecture}}
\newcommand{\econ}{\end{conjecture}}
\newcommand{\bex}{\begin{exercise}}
\newcommand{\eex}{\end{exercise}}
\newcommand{\bcou}{\begin{counterex}}
\newcommand{\ecou}{\end{counterex}}

\newenvironment{Proof}[1][]{\noi\textbf{Proof #1}}{\QED}
\newcommand{\bpro}{\begin{Proof}}
\newcommand{\epro}{\end{Proof}}

\newcommand{\be}{\begin{equation}}
\newcommand{\ee}{\end{equation}}
\newcommand{\ba}{\begin{array}}
\newcommand{\ea}{\end{array}}
\newcommand{\bc}{\be\begin{array}{r@{\,}c@{\,}l}}
\newcommand{\ec}{\end{array}\ee}

\newcommand{\ga}{\gamma}

\newcommand{\de}{\delta}

\newcommand{\eps}{\varepsilon}
\newcommand{\la}{\lambda}
\newcommand{\La}{\Lambda}

\newcommand{\tet}{\theta}

\newcommand{\Ai}{{\cal A}}

\newcommand{\Pc}{{\cal P}}

\newcommand{\Si}{{\cal S}}

\newcommand{\R}{{\mathbb R}}

\newcommand{\Z}{{\mathbb Z}}

\newcommand{\E}{{\mathbb E}}
\renewcommand{\P}{{\mathbb P}}

\newcommand{\up}{\uparrow}
\newcommand{\down}{\downarrow}
\newcommand{\sub}{\subset}
\newcommand{\beh}{\backslash}

\newcommand{\Asto}[1]{\underset{{#1}\to\infty}{\Longrightarrow}}

\newcommand{\ti}{\tilde}
\newcommand{\dgg}{\dagger}
\newcommand{\ov}{\overline}

\newcommand{\ffrac}[2]{{\textstyle\frac{{#1}}{{#2}}}}

\newcommand{\di}{\mathrm{d}}
\newcommand{\half}{{[0,\infty)}}
\newcommand{\expo}{\mbox{\large\it e}}
\newcommand{\ex}[1]{\expo^{\,\textstyle{#1}}}
\newcommand{\ha}{\ffrac{1}{2}}

\newcommand{\Pfp}{\Pc_{{\rm fin},\,+}}
\newcommand{\Pf}{\Pc_{\rm fin}}

\newcommand{\nucirc}{{\accentset{\circ}{\nu}}}

\begin{document}

\makeatletter\@addtoreset{equation}{section}
\makeatother\def\theequation{\thesection.\arabic{equation}} 

\renewcommand{\labelenumi}{{\rm (\roman{enumi})}}
\renewcommand{\theenumi}{\roman{enumi}}

\title{A simple proof of exponential decay\\ of subcritical contact
  processes}
\author{Jan~M.~Swart
\footnote{Institute of Information Theory and
Automation of the ASCR (\' UTIA),
Pod vod\'arenskou v\v e\v z\' i 4,
18208 Praha 8,
Czech Republic;
swart@utia.cas.cz}
}

\date{\today}

\maketitle

\begin{abstract}\noi
This paper gives a new, simple proof of the known fact that for contact
processes on general lattices, in the subcritical regime the expected number
of infected sites decays exponentially fast as time tends to infinity. The
proof also yields an explicit bound on the survival probability below the
critical recovery rate, which shows that the critical exponent associated with
this function is bounded from below by its mean-field value. The main idea of
the proof is that if the expected number of infected sites decays slower than
exponentially, then this implies the existence of a harmonic function that can
be used to show that the process survives for any lower value of the recovery
rate.
\end{abstract}
\vspace{.5cm}

\noi
{\it Acknowledgement.} Work sponsored by grant 16-15238S of the Czech Science Foundation (GA CR).

%


\section{Introduction and results}

A contact process is a Markov process $\eta=(\eta_t)_{t\geq 0}$ taking values in
the subsets of a countable set $\La$, with the following description. If
$i\in\eta_t$, then we say that the site $i$ is infected at time $t$; otherwise
it is healthy. Infected sites $i$ infect healthy sites $j$ with
\emph{infection rate} $a(i,j)\geq 0$, and infected sites become healthy with
\emph{recovery rate} $\de\geq 0$. The formal generator of the process is given
by
\bc\label{Gdef}
Gf(A)&:=&\dis\sum_{i,j\in\La}a(i,j)1_{\{i\in A\}}1_{\{j\notin A\}}
\{f(A\cup\{j\})-f(A)\}\\[5pt]
&&\dis+\de\sum_{i\in\La}1_{\{i\in A\}}\{f(A\beh\{i\})-f(A)\}.
\ec
In the classical set-up, $\La=\Z^d$ and the infection rates are symmetric and
translation-invariant, but other lattices such as regular trees have also
been considered. We refer to \cite{Lig85,Lig99} as a general reference.

In what follows, we will need processes that are translation-invariant in some
sense. A simple way to formalize this, which includes many classical examples
such as processes on $\Z^d$ and regular trees, is to assume that $\La$ is a
group with group action $(i,j)\mapsto ij$, inverse operation $i\mapsto
i^{-1}$, and unit element $0$ (also refered to as the origin). We then assume
that the infection rates $a:\La\times\La\to\half$ satisfy $a(i,i)=0$
$(i\in\La)$ and
\be\ba{rl}\label{assum}
{\rm (i)}&a(i,j)=a(ki,kj)\qquad\qquad(i,j,k\in\La),\\[5pt]
{\rm (ii)}&\dis|a|:=\sum_{i\in \La}a(0,i)<\infty.
\ec
Here (i) says that the infection rates are translation invariant (w.r.t.\ to
the left action of the group on itself), while (ii) guarantees that the
process is well-defined \cite[Thm~I.3.9]{Lig85}. In general, we do not assume
that the infection rates are symmetric, i.e., we allow for the case that
$a\neq a^\dgg$ where we define \emph{reversed infection rates} as
$a^\dgg(i,j):=a(j,i)$.  Using notation as in \cite{Swa09,SS14}, we call the
process with generator in (\ref{Gdef}) the \emph{$(\La,a,\de)$-contact
  process}.

It is well-known \cite[Thm~VI.1.7]{Lig85} that the $(\La,a,\de)$-contact
process $\eta$ and the $(\La,a^\dgg,\de)$-contact process $\eta^\dgg$ are dual
in the sense that
\be\label{dual}
\P[\eta^A_t\cap B\neq\emptyset]=\P[A\cap\eta^{\dgg\,B}_t\neq\emptyset]
\qquad(A,B\sub\La,\ t\geq 0),
\ee
where $\eta^A_t$ and $\eta^{\dgg\,B}_t$ denote the processes started in
$\eta^A_0=A$ and $\eta^{\dgg\,B}_0=B$, respectively.


We say that the $(\La,a,\de)$-contact process \emph{survives} if
$\P[\eta^A_t\neq\emptyset\ \forall t\geq 0]>0$
for some, and hence for all finite nonempty $A$.
We let
\be
\tet(\La,a,\de):=\P\big[\eta^{\{0\}}_t\neq\emptyset\ \forall t\geq 0\big]
\ee
denote the survival probability started from a single infected site,
and call
\be\label{dec}
\de_{\rm c}=\de_{\rm c}(\La,a):=\sup\big\{\de\geq 0:\tet(\La,a,\de)>0\big\}
\ee
the \emph{critical recovery rate}. It is known that $\de_{\rm c}<\infty$. If
$\La$ is finite, then $\de_{\rm c}=0$, but if $\La$ is infinite, then it is
often the case that $\de_{\rm c}>0$. In particular, this is true if $\La$ is
finitely generated and $a$ satisfies a weak irreducibility condition
\cite[Lemma~4.18]{Swa07}. For non-finitely generated infinite groups,
irreducibility is in general not enough to guarantee $\de_{\rm c}>0$
\cite{AS10}. It is well-known that
\be
\P\big[\eta_t^\La\in\cdot\,\big]\Asto{t}\ov\nu,
\ee
where $\ov\nu$ is an invariant law of the $(\La,a,\de)$-contact process, known
as the \emph{upper invariant law}. Using duality, it is not hard to prove that
$\ov\nu=\de_\emptyset$ if the dual $(\La,a^\dgg,\de)$-contact process dies
out, while $\ov\nu$ is concentrated on the nonempty subsets of $\La$ if the
dual process survives \cite[Thms~VI.1.6 and 1.10]{Lig85}. In the latter case,
we say that $\ov\nu$ is \emph{nontrivial}.

It follows from subadditivity (see \cite[Lemma~1.1]{Swa09}) that for any
$(\La,a,\de)$-contact process, there exists a constant $r=r(\La,a,\de)$ with
$-\de\leq r\leq|a|-\de$ such that
\be\label{rdef}
r=\lim_{t\to\infty}\ffrac{1}{t}\log\E\big[|\eta^A_t|\big]\qquad
\mbox{for all finite nonempty }A\sub\La.
\ee
We call $r$ the \emph{exponential growth rate}. The following simple
properties of $r$ are proved in \cite[Theorem~1.2]{Swa09}:
\be\ba{rl}\label{rprop}
{\rm(i)}&\dis r(\La,a,\de)=r(\La,a^\dgg,\de),\\[5pt]
{\rm(ii)}&\dis
\mbox{The function $\de\to r(\La,a,\de)$ is nonincreasing and Lipschitz}\\
&\dis\mbox{continuous on $\half$, with Lipschitz constant 1.}\\[5pt]
\ec

The main aim of the present paper is to present a new, simple proof of the
following known fact.

\bt[Sharpness of the phase transition]\label{T:sharp}
For any $(\La,a,\de)$-contact process, one has $r(\La,a,\de)<0$ if and only if
$\de>\de_{\rm c}$.
\et

We note that Theorem~\ref{T:sharp} and formula (\ref{rprop})~(i) together
imply that $\de_{\rm c}(\La,a)=\de_{\rm c}(\La,a^\dgg)$. By duality, it
follows that for a $(\La,a,\de)$-contact process, the critical points for
survival and nontriviality of the upper invariant law are the same,
which in our present general setting is a nontrivial fact.

The well-known graphical representation of the contact process shows that it
is, in its essence, a form of oriented percolation. Theorem~\ref{T:sharp}
then says that in the whole subcritical regime, connection probabilities decay
exponentially fast.

Historically, such statements were first proved for unoriented percolation, by
Menshikov \cite{Men86} and by Aizenman and Barsky \cite{AB87}; both proofs can
be found in \cite{Gri99}. The proof of \cite{AB87} is based on differential
inequalities involving two parameters: the percolation parameter and the
strength of an external field. Recently, Duminil-Copin and Tassion
\cite{DT1,DT2} have found a much simpler proof which depends on a single
differential inequality and no longer requires the introduction of an external
field.

For contact processes on $\Z^d$, Theorem~\ref{T:sharp} was first proved by
Bezuidenhout and Grimmett \cite{BG91}, who adopted the method of \cite{AB87}
to the oriented, continuous-time setting. This has been generalized to
processes on general transitive graphs in \cite{AJ07}; their arguments also
carry over to general $(\La,a,\de)$-contact processes, as spelled out in the
appendix of \cite{SS14}. The proof of Duminil-Copin and Tassion \cite{DT1,DT2}
works for oriented percolation as well. With a bit of work, it is likely it
can also be adapted to the continuous-time setting of the contact process.

The new proof of Theorem~\ref{T:sharp} presented here is quite different from
the previous proofs. We will see that the assumption that $r(\La,a,\de)=0$
implies the existence of an, in general infinite, invariant measure for the
$(\La,a^\dgg,\de)$-contact process, that gives rise to a harmonic function for
the dual $(\La,a,\de)$-contact process. Lowering the recovery rate a bit turns
this harmonic function into a subharmonic function that allows one to prove
that the $(\La,a,\de-\eps)$-contact process survives for any $\eps>0$. This
method cannot easily be adapted to unoriented percolation, but on the other
hand there seems to be hope that it may be applied to more general interacting
particle systems.

Proofs of sharpness of the phase transition using differential inequalities
typically yield as a side result that the critical exponent associated with
the function $\tet$ is bounded from below by its mean-field value 1; compare,
e.g., \cite[formula~(1.15)]{BG91}, \cite[Part~1 of Thm~1.1]{DT1} or
\cite[Lemma~A.2]{SS14}. Our proof also yields such a result and in fact leads
to the following explicit bound.

\bt[Lower bound on survival probability]\label{T:lwbd}
Let $\phi:(0,1)\to(0,1)$ be implicitly defined by
\be\label{phiga}
\phi(\ga):=1-e^{-\eps}\quad\mbox{where}\quad
\ga=\frac{\eps+e^\eps\eps}{2+e^\eps\eps}
\qquad(0<\eps<2).
\ee
Then $\phi(\ga)=\ga-\frac{1}{2}\ga^2+O(\ga^3)$ as $\ga\to 0$, and
\be\label{lwbd}
\tet\big(\La,a,(1-\ga)\de_{\rm c}\big)\geq\phi(\ga)\qquad(0<\ga<1).
\ee
\et

\section{Proofs}

Let $\Pc=\Pc(\La):=\{A:A\sub\La\}$ denote the set of all subsets of $\La$. We
also set $\Pc_+:=\{A\in\Pc:A\neq\emptyset\}$, $\Pf:=\{A\in\Pc:|A|<\infty\}$,
and $\Pfp:=\Pc_+\cap\Pf$, where $|A|$ denotes the cardinality of a set $A$. We
equip $\Pc\cong\{0,1\}^\La$ with the product topology, making it into a
compact space. Now $\Pc_+$, being a punctured version of $\Pc$, is locally
compact. Recall that a measure on a locally compact space is locally finite if
it gives finite mass to compact sets. By \cite[Lemma~1.1]{SS14}, a measure
$\mu$ on $\Pc_+$ is locally finite if and only if it satisfies one, and hence
both of the following equivalent conditions:
\be
{\rm(i)}\int\mu(\di A)1_{\{i\in A\}}<\infty\ \forall i\in\La\quad
{\rm(ii)}\int\mu(\di A)1_{\{A\cap B\neq\emptyset\}}<\infty\ \forall B\in\Pfp.
\ee

For $A\sub\La$ and $i\in\La$, we write $iA:=\{ij:j\in A\}$, and for any
$\Ai\sub\Pc$ we write $i\Ai:=\{iA:A\in\Ai\}$. We say that a measure $\mu$ on
$\Pc$ is (spatially) \emph{homogeneous} if it is invariant under the left
action of the group, i.e., if $\mu(\Ai)=\mu(i\Ai)$ for each $i\in\La$ and 
measurable $\Ai\sub\Pc$.

It is possible to evolve locally finite (but possibly infinite) starting
measures according to the time evolution of a contact process, as follows. For
a given $(\La,a,\de)$-contact process, we define subprobability kernels $P_t$
$(t\geq 0)$ on $\Pc_+$ by
\be\label{Ptdef}
P_t(A,\,\cdot\,):=\P\big[\eta^A_t\in\cdot\,\big]\big|_{\Pc_+}
\qquad(t\geq 0,\ A\in\Pc_+),
\ee
where $|_{\Pc_+}$ denotes restriction to $\Pc_+$, and we define $P^\dgg_t$
similarly for the dual $(\La,a^\dgg,\de)$-contact process. For any measure
$\mu$ on $\Pc_+$, we write
\be
\label{muPtdef}
\mu P_t:=\int\mu(\di A)P_t(A,\,\cdot\,)\qquad(t\geq 0),
\ee
which is the restriction to $\Pc_+$ of the `law' at time $t$ of the
$(\La,a,\de)$-contact process started in the initial (possibly infinite) `law'
$\mu$. If $\mu$ is a homogeneous, locally finite measure on
$\Pc_+$, then $\mu P_t$ is a homogeneous, locally finite measure on $\Pc_+$
for each $t\geq 0$ (see \cite[Lemma~3.3]{Swa09} or \cite[Lemma~2.4]{SS14}).

Following \cite{Swa09}, we say that a measure $\mu$ on $\Pc_+$ is an
\emph{eigenmeasure} of an $(\La,a,\de)$-contact process if $\mu$ is nonzero,
locally finite, and there exists a constant $\la\in\R$ such that
\be\label{eigen}
\mu P_t=e^{\la t}\mu\qquad(t\geq 0).
\ee
We call $\la$ the associated \emph{eigenvalue}. We cite the following result
from \cite[Prop.~1.4]{Swa09}.

\bl[Existence of eigenmeasures]\label{L:eigex}
Each $(\La,a,\de)$-contact process has a homogene\-ous eigenmeasure $\nucirc$
with eigenvalue $r=r(\La,a,\de)$.
\el

This is proved in \cite{Swa09} along the following lines. First, it is shown
that the function $t\mapsto\E[|\eta^{\{0\}}_t|]$ is submultiplicative, which
by Fekete's lemma implies that there exists a constant $r$ such that
\be\label{Fekete}
\E[|\eta^{\{0\}}_t|]=e^{rt+o(t)},
\ee
where $o(t)\geq 0$ and $\lim_{t\to\infty}o(t)/t=0$. Next, one defines
\be
\nu_\la:=\int_0^\infty\!\mu P_t\,e^{-\la t}\di t
\quad(\la>r)
\quad\mbox{with}\quad
\mu:=\sum_{i\in\La}\de_{\{i\}}.
\ee
Then $\nu_\la G=\la\nu_\la-\mu$ and $\int\nu_\la(\di A)1_{\{i\in
  A\}}=\int_0^\infty\E[|\eta^{\{0\}}_t|]e^{-\la t}\di t=:\pi_\la$. Since the
function $o(t)$ in (\ref{Fekete}) is nonnegative, we have that
$\pi_\la\up\infty$ as $\la\down r$. Using this, it can be shown that the
normalized measures $\pi_\la^{-1}\nu_\la$ converge vaguely along some
subsequence $\la_n\down r$ to a limit $\nucirc$ which satisfies $\nucirc
G=r\nucirc$, i.e., $\nucirc$ is an eigenmeasure with eigenvalue $r$.

In \cite{Swa09,SS14}, it is moreover proved that the eigenmeasure $\nucirc$ from
Lemma~\ref{L:eigex} is in certain situations unique up to a multiplicative
constant. For our present purposes, however, we only need existence.


It is well-known that if two Markov processes are dual, then invariant laws of
one Markov process give rise to harmonic functions for its dual. A similar
statement holds for eigenmeasures. We cite the following lemma from
\cite[Lemma~3.5]{Swa09}.

\bl[Eigenfunctions]\label{L:eigfun}
Let $\mu$ be a homogeneous eigenmeasure of a $(\La,a,\de)$-contact process with
eigenvalue $\la$. Then
\be\label{eigfun}
h(A):=\int\mu(\di B)1_{\txt\{A\cap B\neq\emptyset\}}\qquad(A\in\Pf)
\ee
satisfies $G^\dgg h=\la h$, where $G^\dgg$ is defined as in
(\ref{Gdef}) but for the dual $(\La,a^\dgg,\de)$-contact process.
\el

\bpro[Proof of Theorem~\ref{T:sharp}]
Let $\de'_{\rm c}:=\inf\{\de\geq 0:r(\La,a,\de)<0\}$. By (\ref{rprop})~(ii)
and the bounds $-\de\leq r\leq|a|-\de$, we have $r(\La,a,\de'_{\rm c})=0$ and
$r(\La,a,\de)<0$ for all $\de>\de'_{\rm c}$. It is easy to see that the
$(\La,a,\de)$-contact process dies out when $r<0$, so to complete the proof
it suffices to prove that the $(\La,a,\de)$-contact process
survives for $\de<\de'_{\rm c}$.

By (\ref{rprop})~(i), $r(\La,a^\dgg,\de'_{\rm c})=0$, so Lemma~\ref{L:eigex}
tells us that the $(\La,a^\dgg,\de'_{\rm c})$-contact process has a
homogeneous eigenmeasure $\nucirc^\dgg$ with eigenvalue $0$. Since
$\nucirc^\dgg$ is nonzero, concentrated on $\Pc_+$, and homogeneous, we can
normalize $\nucirc^\dgg$ such that $\int\nucirc^\dgg(\di B)1_{\{i\in B\}}=1$
$(i\in\La)$. With this normalization, we define $h:\Pf\to\R$ by
\be\label{hdef}
h(A):=\int\nucirc^\dgg(\di B)1_{\txt\{A\cap B\neq\emptyset\}}\qquad(A\in\Pf).
\ee
Letting $G$ denote the generator of the $(\La,a,\de'_{\rm c})$-contact
process, Lemma~\ref{L:eigfun} tells us that $Gh=0$, i.e., $h$ is a harmonic
function.

Let $\eps_1,\eps_2>0$ be constants to be chosen later and let $\ti G$ denote
the generator of the $(\La,(1+\eps_1)a,(1-\eps_2)\de'_{\rm c})$-contact
process. We use Lemma~\ref{L:posdrift} in the appendix to transform $h$ into a
bounded function $f_\eps$ as in (\ref{feps}) such that $f_\eps$ is subharmonic
for $\ti G$. To this aim, we check condition (\ref{posdrift}). Using the fact
that $Gh=0$, we see that
\bc
\ti Gh(A)&=&\dis\eps_1\sum_{i,j\in\La}a(i,j)1_{\{i\in A\}}1_{\{j\notin A\}}
\big(h(A\cup\{j\})-h(A)\big)\\[5pt]
&&\dis+\eps_2\de'_{\rm c}\sum_{i\in\La}1_{\{i\in A\}}\big(h(A)-h(A\beh\{i\})\big).
\ec
Let $\ti H_\eps$ be defined as in (\ref{Hphi}) but with $\ti G$ instead of
$G$, i.e.,
\bc
\ti H_\eps h(A)
&=&\dis(1+\eps_1)\sum_{i,j\in\La}a(i,j)1_{\{i\in A\}}1_{\{j\notin A\}}
\phi_\eps\big(h(A\cup\{j\})-h(A)\big)\\[5pt]
&&\dis+(1-\eps_2)\de'_{\rm c}\sum_{i\in\La}1_{\{i\in A\}}
\phi_\eps\big(h(A\beh\{i\})-h(A)\big).
\ec
It follows from (\ref{hdef}) that $h$ is an increasing function in the sense
that $A\sub A'$ implies $h(A)\leq h(A')$. Moreover, combining (\ref{hdef})
with our normalization $\int\nucirc^\dgg(\di B)1_{\{i\in B\}}=1$ we see that
\be
h(A)-h\big(A\beh\{i\}\big)
=\int\nucirc^\dgg(\di B)1_{\txt\{A\cap B=\{i\}\}}\leq 1\qquad(i\in A\in\Pfp).
\ee
Since $\phi_\eps(0)=0$, $\phi'_\eps(0)=0$, and $\phi''_\eps(z)=\eps e^{-\eps
  z}$, we see that
\be
\phi_\eps(z)\leq\ha\eps z^2\quad(0\leq z\leq 1)
\quand
\phi_\eps(z)\leq\ha\eps e^\eps z^2\quad(-1\leq z\leq 0).
\ee
In view of this, using moreover that $z^2\leq|z|$ for $|z|\leq 1$,
we can estimate
\be\ba{l}
\dis\big(\ti Gh-\ti H_\eps h\big)(A)\\[5pt]
\dis\quad\geq\big[\eps_1-\ha\eps(1+\eps_1)\big]
\sum_{i,j\in\La}a(i,j)1_{\{i\in A\}}1_{\{j\notin A\}}\big(h(A\cup\{j\})-h(A)\big)\\[5pt]
\dis\quad\phantom{\leq}
+\big[\eps_2-\ha\eps e^\eps(1-\eps_2)\big]
\de'_{\rm c}\sum_{i\in\La}1_{\{i\in A\}}\big(h(A)-h(A\beh\{i\})\big),
\ec
so condition (\ref{posdrift}) is satisfied and hence $\ti Gf_\eps\geq 0$ when
we choose $\eps_1,\eps_2$ in such a way that
\be\label{epschoice}
\frac{\eps_1}{1+\eps_1}=\ha\eps
\quand
\frac{\eps_2}{1-\eps_2}=\ha\eps e^\eps.
\ee
Let $(\ti\eta_t)_{t\geq 0}$ denote the process with generator $\ti G$, started in
$\ti\eta_0=\{0\}$, i.e., with a single infected site at the origin. It is easy
to see that if a contact process with infection rates satisfying (\ref{assum})
is started in a finite initial state, then it stays finite for all time, so 
$(\ti\eta_t)_{t\geq 0}$ is nonexplosive. Since $\ti Gf_\eps\geq 0$ and
since $f_\eps$ is a bounded function, we have that $f_\eps(\ti\eta_t)$
is a bounded submartingale that converges to an a.s.\ limit
$\lim_{t\to\infty}f_\eps(\ti\eta_t)=:F_\infty$ with
\be
\E[F_\infty]\geq f_\eps(\{0\})=\frac{1}{\eps}(1-e^{-\eps}),
\ee
where we have used that $h(\{0\})=1$ by our normalization of $\nucirc^\dgg$.
Since $f_\eps(\emptyset)=0$ and $f_\eps\leq\eps^{-1}$, we have
\be\label{survbd}
\P[\ti\eta_t\neq\emptyset\ \forall t\geq 0]\geq\P[F_\infty>0]\geq\eps\E[F_\infty]
\geq 1-e^{-\eps}.
\ee
By a trivial rescaling of time, it follows that
\be\label{tetlow}
\tet\big(\La,a,\ffrac{1-\eps_2}{1+\eps_1}\de'_{\rm c}\big)
=\tet\big(\La,(1+\eps_1)a,(1-\eps_2)\de'_{\rm c}\big)\geq 1-e^{-\eps},
\ee
where $\eps_1,\eps_2$ are defined in terms of $\eps$ as in (\ref{epschoice}).
Since $\eps>0$ is arbitrary, this completes the proof. In
particular, our argument shows that $\de'_{\rm c}=\de_{\rm c}$.
\epro

\bpro[Proof of Theorem~\ref{T:lwbd}]
We observe that by (\ref{epschoice})
\be
\frac{1}{1+\eps_1}=1-\ha\eps
\quand
\frac{1}{1-\eps_2}=1+\ha\eps e^\eps.
\ee
Defining $\ga$ as in (\ref{phiga}), we have that
\be
\frac{1-\eps_2}{1+\eps_1}=\frac{1-\ha\eps}{1+\ha\eps e^\eps}=1-\ga.
\ee
Then (\ref{tetlow}) says that
\be\label{tetlow2}
\tet\big(\La,a,(1-\ga)\de_{\rm c}\big)\geq 1-e^{-\eps},
\ee
which is (\ref{lwbd}). Since
\be
\ga=\eps+O(\eps^3)
\quand
1-e^{-\eps}=\eps-\ha\eps^2+O(\eps^3)
\quad\mbox{as }\eps\to 0,
\ee
we see that $\phi(\ga)=\ga-\ha\ga^2+O(\ga^3)$ as $\ga\to 0$.
\epro

\appendix

\section{Transformation of submartingales}

Let $\Si$ be a countable set and let $G$ be a so-called Q-matrix on $\Si$,
i.e., $(G(x,y))_{x,y\in\Si}$ are real constants such that $G(x,y)\geq 0$ for
$x\neq y$ and $\sum_{y\in\Si}G(x,y)=0$. For any real function $f$ on
$\Si$, we write
\be
Gf(x):=\sum_{y\in\Si}G(x,y)f(y)
=\sum_{y\in\Si}G(x,y)\big(f(y)-f(x)\big)
\qquad(x\in\Si),
\ee
whenever the infinite sums are well-defined. Then $G$ is the the generator of
a (possibly explosive) continuous-time Markov chain $(X_t)_{t\geq 0}$ in
$\Si$. A function $h$ such that $Gh\geq 0$ is called \emph{subharmonic}. The
following simple lemma says, roughly speaking, that an unbounded, nonnegative
subharmonic function that has a sufficiently positive drift and not too
large fluctuations can be transformed into a bounded subharmonic function.


\bl[Transformation of submartingales]\label{L:posdrift}
Let $h$ be a real function on $\Si$ and let $\eps>0$. Then the function
\be\label{feps}
f_\eps(x):=\frac{1}{\eps}\big(1-\ex{-\eps h(x)}\big)\qquad(x\in\Si)
\ee
satisfies $Gf_\eps\geq 0$ if and only if
\be\label{posdrift}
Gh-H_\eps h\geq 0,
\ee
where
\be\label{Hphi}
H_\eps h(x):=\sum_{y\in\Si}G(x,y)\phi_\eps\big(h(y)-h(x)\big)
\quad\mbox{with}\quad
\phi_\eps(z):=\eps^{-1}(e^{-\eps z}-1+\eps z).
\ee
\el 
\bpro
Let $g_\eps(z):=\eps^{-1}(1-e^{-\eps z})$ $(z\in\R)$. Then, for any $z,z_0\in\R$,
\be
g_\eps(z)=g_\eps(z_0)+\big\{(z-z_0)-\phi_\eps(z-z_0)\big\}e^{-\eps z_0}.
\ee
It follows that
\bc
\dis Gf_\eps(x)&=&\dis\sum_{y\in\Si}G(x,y)
\big\{g_\eps\big(h(y)\big)-g_\eps\big(h(x)\big)\big\}\\[5pt]
&=&\dis\ex{-\eps h(x)}\sum_{y\in\Si}G(x,y)
\big\{\big(h(y)-h(x)\big)-\phi_\eps\big(h(y)-h(x)\big)\big\},
\ec
which is nonnegative if and only if (\ref{posdrift}) holds.
\epro


\begin{thebibliography}{DT15b}

\bibitem[AB87]{AB87}
M.~Aizenman and D.J.~Barsky.
Sharpness of the phase transition in percolation models.
\emph{Comm. Math. Phys.}~{108}, (1987), 489--526.

\bibitem[AJ07]{AJ07}
M.~Aizenman and P.~Jung.
On the critical behavior at the lower phase transition of the contact
process.
\emph{Alea}~{3}, (2007), 301--320.

\bibitem[AS10]{AS10}
S.R.~Athreya and J.M.~Swart.
Survival of contact processes on the hierarchical group.
\emph{Prob.\ Theory Relat.\ Fields}~{147(3)}, (2010), 529-563.

\bibitem[BG91]{BG91}
C.~Bezuidenhout and G.~Grimmett.
Exponential decay for subcritical contact and percolation processes.
\emph{Ann.\ Probab.}~{19(3)}, (1991), 984--1009.  

\bibitem[DT15a]{DT1}
H.~Duminil-Copin and V.~Tassion.
A new proof of the sharpness of the phase transition for Bernoulli percolation
and the Ising model.
\emph{Commun.\ Math.\ Phys.}~343(2), (2016), 725--745, 

\bibitem[DT15b]{DT2}
H.~Duminil-Copin and V.~Tassion.
A new proof of the sharpness of the phase transition for Bernoulli percolation
on $\Z^d$.
Preprint (2015), ArXiv:1502.03051.
To appear in \emph{L'Enseignement Math\'ematique}.

\bibitem[Gri99]{Gri99}
G.~Grimmett.
\emph{Percolation 2nd ed.},
Vol.~{321} Grundlehren der Mathematischen Wissenschaften.
Springer-Verlag, Berlin, 1999.

\bibitem[Lig85]{Lig85}
T.M.~Liggett.
\emph{Interacting Particle Systems.}
Springer-Verlag, New York, 1985.

\bibitem[Lig99]{Lig99}
T.M.~Liggett.
\emph{Stochastic interacting systems: contact, voter and exclusion processes.}
Springer-Verlag, Berlin, 1999.

\bibitem[Men86]{Men86}
M.V.~Menshikov.
Coincidence of the critical points in percolation problems.
\emph{Soviet Math. Dokl.}~{33}, (1986), 856--859.

\bibitem[SS14]{SS14}
A.~Sturm and J.M.~Swart.
Subcritical contact processes seen from a typical infected site.
\emph{Electron.\ J.\ Probab.}~{19} (2014), no.~53, 1--46.

\bibitem[Swa07]{Swa07}
J.M.~Swart.
Extinction versus unbounded growth.
Habilitation Thesis of the University Erlangen-N\"urnberg, 2007.
ArXiv:math/0702095v1.

\bibitem[Swa09]{Swa09}
J.M.~Swart.
The contact process seen from a typical infected site.
\emph{J.\ Theoret.\ Probab.}~{22(3)}, (2009), 711-740.

\end{thebibliography}
\end{document}